\documentclass[11pt]{amsart}
\usepackage{amsmath,amssymb,amsbsy,amsfonts,amsthm,latexsym,amsopn,amstext,amsxtra,epic, euscript,amscd,indentfirst}
\begin{document}

\newtheorem{theorem}{Theorem}
\newtheorem{conjecture}[theorem]{Conjecture}
\newtheorem{proposition}[theorem]{Proposition}
\newtheorem{question}[theorem]{Question}
\newtheorem{lemma}[theorem]{Lemma}
\newtheorem{cor}[theorem]{Corollary}
\newtheorem{obs}[theorem]{Observation}
\newtheorem{proc}[theorem]{Procedure}
%% DEFINITIONS
\def\Z{\mathbb Z}
\def\Fq{{\mathbb F}_q}
\def\R{\mathbb R}
\def\N{\mathbb N}

\title{Revisiting Toom's proof \\ of Bulgarian Solitaire}

\author{Therese A. Hart}
\address{TAH: Department of Mathematics and Computer Science, Eastern Connecticut State University, Willimantic, CT 06226}
\email{hartt@easternct.edu}
\author{Gabriel J. H. Khan}
\address{GJHK: Boston University, Boston, MA 02215}
\email{GK:gkhan@bu.edu}
\author{Mizan R. Khan}
\address{MRK: Department of Mathematics and Computer Science, Eastern Connecticut State University, Willimantic, CT 06226}
\email{khanm@easternct.edu}

\date{\today}

\begin{abstract}
In this article we give an exposition of Toom's proof of Bulgarian Solitaire that appeared in \emph{Kvant}. We provide more details. We also show how an application of the Chinese Remainder Theorem allows us to generalize the proof. 
\end{abstract}

\maketitle

\section{Introduction}

The following \emph{literary} version of Bulgarian Solitaire appeared as a problem (numbered M655) in the Russian high school journal \emph{Kvant}~\cite{K1} in 
1980. Its solution by Andrei Toom was published in the same journal~\cite{K2} in 1981. We refer the interested reader who knows Russian to the additional source~\cite[144-147]{VGRT}, where the same solution appears.

\vspace{0.1in}

M655. On the table of a clerk at the Circumlocution Office there are  $n$  volumes of Encyclopedia Britannica, ordered in several piles.  Every day, when arriving at work, the clerk takes one volume from each pile and puts all of these in a new pile, after which he rearranges the piles according to the number of volumes in each (in non-increasing order), and fills out a form, recording the number of volumes in each pile.  He never does anything else, except this.
\begin{enumerate}
\item What record will enter the form after a month, if the total number of volumes is $n=3$, $n=6$, $n=10$ (the initial distributions of volumes is arbitrary)?
\item Prove that if the total number of volumes is  $n=k(k+1)/2$, where $k$ is a natural number, then after a certain amount of days the form will start filling up with identical records.
\item* Investigate what will happen after many working days for other values of  $n$.
\end{enumerate} 
\hspace{3.5in} S. Limanov, A. Toom

\vspace{0.2in}

We note that the idea of the \emph{Circumlocution Office} is from the Dickensian novel \emph{Little Dorrit}~\cite[Chapter 10, Book 1]{D}. Martin Gardner popularized Bulgarian Solitaire in the West in his Scientific American column~\cite{MG}. We quote the following from~\cite[page 18]{MG}: ``Our last example of a task that ends suddenly in a counterintuitive way is one you will enjoy modeling with deck of cards. Its origin is unknown, but Graham, who told me about it, says that European mathematicians call it Bulgarian solitaire for reasons he has not been able to discover." Later in the column Gardner credits Jorgen Brandt\cite{JB} for the first proof of the solution. However, the proof by Toom predates this, and more importantly is much simpler. In this article we describe Toom's proof, but include more detail, and introduce a generalization. 

A very nice discussion of the the origins of Bulgarian Solitaire and subsequent literature is given in the preprint~\cite{H}. In particular, the mystery behind its name is solved there. 

We would like to thank both Z. Pozdynakova and C. Yankov for their superb translating services, especially since none of the three authors have the faintest idea of Russian.

\section{Definitions, Notation, Pictures and Statement of Results}

We recast the problem of Bulgarian solitaire into the language of partitions and dynamical systems. Consequently, we begin by stating  some basic definitions and introducing some notation. In this context we follow the usage in~\cite{S}.

A partition $\pi$ of a positive integer $n$ is a nonincreasing sequence  of nonnegative integers such that the sum of terms is $n$, that is,
$$ \pi=\left(\pi_1,\pi_2,\pi_3,\ldots\right) $$
where $ \pi_i $ is a nonnegative integer, $\pi_i \geq \pi_{i+1}$ for all $i$, and $\sum_{i=1}^\infty \pi_i = n.$

Since $\pi_i$ is a nonnegative integer, only a finite number of terms of $\pi$ are non-zero. These are called the parts of $\pi$, and the number of
parts of $\pi$ is called the length of $\pi$ and is denoted $l(\pi)$. In writing specific partitions one typically omits the infinite tail of 
zeros. For example, the partition $(3,3,2,2,1,0,0,\ldots)$ is written as $(3,3,2,2,1)$.

We denote the set of all partitions by $\mathcal P$, and the subset of the partitions of $n$ by ${\mathcal P}_n$. For $\lambda,\pi \in {\mathcal P}$ we define the union $\lambda \cup \pi$ to be the partition obtained by merging the entries of $\lambda$ with those of $\pi$ and arranging the resulting entries in nonincreasing order, for example,
$$(4,3,3,2,1) \cup (5,4,3,1,1) = (5,4,4,3,3,3,2,1,1,1).$$

We now define a map $ T: {\mathcal P} \rightarrow {\mathcal P}$ via
\begin{equation}
\label{eq:defn-T}
T(\pi) = (\pi_1 -1, \pi_2 -1, \ldots, \pi_{l(\pi)} -1,0,0,\ldots) \cup (l(\pi),0,0,\dots)
\end{equation}
We note that if $\pi \in {\mathcal P}_n$ then $T(\pi) \in {\mathcal P}_n$. In this paper we study the dynamics of this map $T$. If we start with a partition $\pi \in {\mathcal P}_n$ and look at the iterates of $T$ we get the orbit 
\begin{equation*} {\mathcal O}_T(\pi) = \{ \pi, T(\pi), T^2(\pi), T^3(\pi), \ldots \}. 
\end{equation*}
For fixed value of $n$, the number of partitions of $n$ is finite. Consequently the orbit, ${\mathcal O}_T(\pi)$, is closed and  ends in a $N$-cycle. The problem of Bulgarian solitaire is characterizing the periodic points of $T$.

We visualize a partition in a couple of different ways. The first as a Ferrers graph. A Ferrers graph represents a partition as a pattern of dots with the $k$-th row having the same number of dots as the $k$-th term of the partition. So for example, the Ferrers graph of $(5,4,1,1,1)$ is 
\begin{equation*}
\begin{array}{ccccc}
\bullet & \bullet & \bullet & \bullet & \bullet \\
\bullet & \bullet & \bullet & \bullet &  \\
\bullet &  &  &  &  \\
\bullet &  &  &  &  \\
\bullet &  &  &  &  
\end{array}
\end{equation*}

We denote the Ferrers graph of $\pi$ by $F_\pi$. The solution to Bulgarian solitaire involves viewing $F_\pi$ as an arrangement of checkers in the upper triangle of an infinite checkerboard. For example, the partition $(5,4,1,1,1)$ has the following arrangement.
\begin{equation*}
\begin{array}{|c|c|c|c|c|c|c}
\hline 
\bullet & \bullet & \bullet & \bullet & \bullet & \  & \  \\
\hline 
\bullet & \bullet & \bullet & \bullet &   &   & \\
\hline
\bullet &  &  &  &  &  &  \\
\hline
\bullet &  &  &  &  &  &  \\
\hline
\bullet &  &  &  &  &  \\
\hline
&  &  &  &  &  & \\
\hline
&  &  &  &  &  & \\
\end{array}
\end{equation*}

We use the ordered pair $(i,j)$ to denote the square on the $i$-th row and $j$-th column of the checkerboard, and we use $diag[k]$ to denote the $k$-th upper diagonal of the infinite checkerboard, that is,  $diag[k]$ are the squares 
$$diag[k] = \{ \, (1,k),(2,k-1),(3,k-2), \ldots, (k-1,2), (k,1) \}. $$ 
For example, in the picture below, $diag[3]$ is the set of three squares containing the checkers. 
\begin{equation*}
\begin{array}{|c|c|c|c|c}
\hline 
\ & \ & \bullet & \ & \   \\
\hline 
\  & \bullet & \ & \ & \\
\hline
\bullet & \ & \ & \ &    \\
\hline
\ & \ & \ &  \ &   \\
\hline
\ & \ & \ & \ & \\
\end{array}
\end{equation*}
\textbf{Note:} For the remainder of this paper we use the word \emph{diagonal} in a very restrictive sense. It will be only used to reference the sets $diag[k]$. 

\vspace{.1in}

\noindent We now state the solution to Bulgarian solitaire.

\begin{theorem}[Bulgarian Solitaire]
\label{Main-Result}
Let $\pi$ be a partition of $n$, and let $m$ be the integer such that $m(m+1)/2 \le n < (m+1)(m+2)/2$.  $\pi$ is a periodic point of $T$ if and only if 
\begin{equation}\label{eq:per cond}
\# \textrm{ of checkers on } diag[k] = \left\{ \begin{array}{ll} k, & k=1,2, \ldots, m \\
n-m(m+1)/2, & k = m+1 \\
0, & k > m+1. \end{array}\right.
\end{equation} 
So in the special case of $n=m(m+1)/2$ we have that, with respect to $T$, ${\mathcal P}_n$ contains only one periodic point, namely, 
$$(m, m-1, m-2, \dots, 2,1).$$
\end{theorem}

Below is the Ferrer's graph of (5,4,4,2,2,1). The above theorem says that this is a periodic point of $T$.
\begin{equation*}
\begin{array}{ccccc}
\bullet & \bullet & \bullet & \bullet & \bullet \\
\bullet & \bullet & \bullet & \bullet &  \\
\bullet & \bullet & \bullet  & \bullet &  \\
\bullet & \bullet &  &  &  \\
\bullet & \bullet &  &  & \\
\bullet & & & &
\end{array}
\end{equation*}

\section{Reinterpreting the action of $T$}

The reason we like to view $F_\pi$ as an arrangement of checkers on a checkerboard is that it introduces a co-ordinate system. This permits a convenient way to decompose the action of $T$ into different components which then leads to a useful visualization of the action of $T$. 

\begin{obs}[Moving Checkers]
Given a partition $\pi$ we arrive at $T(\pi)$ in the following way.
\begin{enumerate}
\item We remove all of the checkers in the first column of the checkerboard.
\item We translate the remaining checkers down one square and then one square to the left.
\item We now place the checkers that we had removed from the first column into the first row. In doing so we take care to place the checker that originally been on the $(k,1)$-th square onto the $(1,k)$-th square.
\item Finally if $l(\pi) < \pi_1-1$, then we move the checkers in the columns numbering $l(\pi) +1$ through $\pi_1-1$ up one square. 
\end{enumerate}
\end{obs}

We illustrate the above procedure using the partition $(5,4,1)$. 
\begin{equation*}
\begin{array}{|c|c|c|c|c|c|c}
\hline 
\bullet & \bullet & \bullet & \bullet & \bullet & \  & \  \\
\hline 
\bullet & \bullet & \bullet & \bullet &   &   & \\
\hline
\bullet &  &  &  &  &  &  \\
\hline
 &  &  &  &  &  &  \\
\hline
&  &  &  &  &  & \\
\end{array}
\end{equation*}
Our first step gives us the pattern 
\begin{equation*}
\begin{array}{|c|c|c|c|c|c|c}
\hline 
\ & \bullet & \bullet & \bullet & \bullet & \  & \  \\
\hline 
 & \bullet & \bullet & \bullet &   &   & \\
\hline
 &  &  &  &  &  &  \\
\hline
&  &  &  &  &  & \\
\end{array}
\end{equation*}
The second step is 
\begin{equation*}
\begin{array}{|c|c|c|c|c|c}
\hline 
 \ & \ & \ & \ & \  & \  \\
\hline 
\bullet & \bullet & \bullet & \bullet & \  & \  \\
\hline 
\bullet & \bullet & \bullet &   &   & \\
\hline
  &  &  &  &  &  \\
\hline
 &  &  &  &  & \\
\end{array}
\end{equation*}
Our third step is 
\begin{equation*}
\begin{array}{|c|c|c|c|c|c}
\hline 
 \bullet & \bullet & \bullet & \ & \  & \  \\
\hline 
\bullet & \bullet & \bullet & \bullet & \  & \  \\
\hline 
\bullet & \bullet & \bullet &   &   & \\
\hline
  &  &  &  &  &  \\
\hline
 &  &  &  &  & \\
\end{array}
\end{equation*}
Our final step is 
\begin{equation*}
\begin{array}{|c|c|c|c|c|c}
\hline 
 \bullet & \bullet & \bullet & \bullet & \  & \  \\
\hline 
\bullet & \bullet & \bullet & \ & \  & \  \\
\hline 
\bullet & \bullet & \bullet &   &   & \\
\hline
  &  &  &  &  &  \\
\hline
 &  &  &  &  & \\
\end{array}
\end{equation*}

\begin{obs}
By describing the action of $T$ in this way we infer the following. 
\begin{itemize}
\item A checker in the $(k,1)$-th square moves to the $(1,k)$-th square. 
\item A checker in the $(i,j)$-th square, with $ i-1 \leq l(\pi)$ and $1 < j $, moves to the $(i+1,j-1)$-th square. 
\item Finally a checker on the 
$(i,j)$-th square, with $ i-1 > l(\pi)$ and $1 < j $, moves to the $(i,j-1)$-th square.
\end{itemize}
\end{obs}

Consequently, a checker on the $k$-th diagonal may move to the $(k-1)$-th diagonal, but it cannot move to the $(k+1)$-th 
diagonal. With these remarks we now arrive at the following description of the action of $T$. 

\begin{obs}[Action of $T$]\label{Act-T} We can view the action of $T$ as moving checkers down the diagonals of the checkerboard. When a checker reaches the square $(k,1)$, then on the next move it jumps up to the $(1,k)$ square. When a checker moves down one square on the diagonal and then if there is an empty square above it, it moves up vertically one square and thus moves into smaller diagonal.
\end{obs}

For this problem we found it convenient visually to rotate the Ferrers graph clockwise by $45^\circ$, with the center of rotation being the top left hand point of the graph. So for example we  view $(5,4,1)$ as

\begin{equation*}
\begin{array}{ccccccccccc}
& & & & & \bullet &  &  &  &  \\
& & & & \bullet & & \bullet & & & \\
& & & \bullet & & \bullet & & \bullet &  &  \\
& & \circ & & \circ & & \bullet & & \bullet &  \\
& \circ & & \circ & & \circ & & \bullet & & \bullet  \\ 
\end{array}
\end{equation*}

\begin{obs}
Each row of the above triangular array represents a diagonal on the checkerboard. We use circles to represent empty squares on the diagonal. When working with such a triangular array we will refer to the $(i,j)$ square on the checkerboard as the $[i+j-1,j]$ cell of the triangular array, that is, the first co-ordinate of the cell refers to the diagonal on which the square belongs.  We use the word \emph{square} when we visualize the checkers lying on a checkerboard and we use the word \emph{cell} when we visualize the checkers lying in a triangular array. We view the action of $T$ as first moving each checker one cell to the left (with the condition that a checker in cell $[s,1]$ moves into the cell $[s,s]$, that is, it loops to the other side), and then possibly moving a checker diagonally up one cell to the right. 
\end{obs}

Let us illustrate with $(5,4,1)$ and show how we get that 
$T((5,4,1))= (4,3,3)$. We start with
\begin{equation*}
\begin{array}{ccccccccccc}
& & & & & \bullet &  &  & &  \\
& & & & \bullet & & \bullet & & & \\
& & & \bullet & & \bullet & & \bullet &  &  \\
& & \circ & & \circ & & \bullet & & \bullet &  \\
& \circ & & \circ & & \circ & & \bullet & & \bullet  \\ 
\end{array}
\end{equation*}
We now move the checkers one cell to the left to obtain
\begin{equation*}
\begin{array}{ccccccccccc}
& & & & & \bullet &  &  & &  \\
& & & & \bullet & & \bullet & & & \\
& & & \bullet & & \bullet & & \bullet &  &  \\
& & \circ & & \bullet & & \bullet & & \circ &  \\
& \circ & & \circ & & \bullet & & \bullet & & \circ  \\ 
\end{array}
\end{equation*}
and then move the checker in cell $[5,4]$ to cell $[4,4]$ to obtain 
\begin{equation*}
\begin{array}{ccccccccccc}
& & & & & \bullet &  &  & &  \\
& & & & \bullet & & \bullet & & & \\
& & & \bullet & & \bullet & & \bullet &  &  \\
& & \circ & & \bullet & & \bullet & & \bullet &  \\
& \circ & & \circ & & \bullet & & \circ & & \circ  \\ 
\end{array}
\end{equation*}
This is the rotated Ferrers graph of $(4,3,3)$. On occasion we adopt the following shortcut in presenting the rotated Ferrers graph by omitting the filled rows. So for example we permit the following to represent $(4,3,3)$.
\begin{equation*}
\begin{array}{ccccccccccc}
& & \circ & & \bullet & & \bullet & & \bullet &  \\
& \circ & & \circ & & \bullet & & \circ & & \circ  \\ 
\end{array}
\end{equation*}
\section{Proof of Theorem~\ref{Main-Result}}

We begin with a small lemma that we will require in our proof. It is an immediate consequence of the Chinese Remainder Theorem.

\begin{lemma}\label{CRTApp}
Let $a,b,m,n,u,v \in \Z$ with $\gcd(m,n)=\gcd(u,m)=\gcd(v,n)=1$. Then for any $c \in \Z$, there exist a corresponding $k \in \Z$ such that
\begin{equation}
c \equiv (a + ku) \!\!\!\! \pmod{m} \textrm{ and } c \equiv (b+kv) \!\!\!\! \pmod{n}.
\end{equation}
\end{lemma}

\begin{proof}
Let $u^\prime,v^\prime \in \Z$ such that $uu^\prime \equiv 1 \pmod{m}$ and $vv^\prime \equiv 1 \pmod{n}$. Consider the simultaneous congruences
$$ x \equiv (c-a)u^\prime \!\!\!\! \pmod{m} \textrm{ and } x \equiv (c-b)v^\prime \!\!\!\! \pmod{n}.$$
By the Chinese Remainder Theorem, this system of simultaneous congruences has a solution $k$.
\end{proof}

\begin{obs}[Key idea of proof] The key idea that we exploit for the proof of Theorem~\ref{Main-Result} is that a checker can never move from a smaller diagonal to a larger diagonal. The reader may first want to check the action of $T$ on $(4,3,3)$. We found this to be an instructive example illustrating our proof.
\end{obs}

\begin{proof}[Proof of Theorem~\ref{Main-Result}] $(\Rightarrow)$ We will prove this direction by proving the contrapositive. %We denote the rotated Ferrers graph of 
%$\pi$ as $F^{\vartriangle}_\pi$. 
We begin by observing that if $diag[d]$ contains an empty square then so does every higher diagonal. Let $\pi \in {\mathcal P}_n$ be a 
partition that does not satisfy ~\eqref{eq:per cond}. Consequently there exists an integer $l$ such that $diag[l]$ and $diag[l+1]$ (of $F_\pi$) contains both checkers and spaces. Without loss of generality we can assume that $l$ is minimal, that is, for every $d<l$ the 
diagonal $diag[d]$ of $F_\pi$ contains $d$ checkers and no spaces.

Let cell $[l,a]$ be empty and let cell $[l+1,b]$ contain a checker. By Lemma~\ref{CRTApp} for each $c\in \Z$ with $1 \leq c \leq  l$, there corresponds a smallest positive integer $k(c)$ such that 
$$ a-k(c) \equiv c\!\!\!\! \pmod{l}, \, b-k(c) \equiv c \!\!\!\! \pmod{(l+1)}. $$

Let $K(C) = \min(\{k(c) \, : \, c=1,2, \ldots, l \}).$ We claim that a checker will move from $diag[l+1]$ to $diag[l]$ after \emph{at most} $K(C)$ iterations. To see this we assume that we have iterated $T$ $(K(C)-1)$ times and at no point a checker has moved  from $diag[l+1]$ to $diag[l]$. Now when we take the $K(C)$-th iteration we find that after translating the checkers to the left (in the triangular array), there will be a checker in the $[l+1,C]$ cell and an empty space in the $[l,C]$ cell. At this point the checker in the $[l+1,C]$ cell will move to the $[l,C]$ cell.

As no checker will move from  $diag[l]$ to $diag[l+1]$, no amount of iterations of $T$ will yield the original partition, and consequently  $\pi$ is not a periodic point of $T$.
\vspace{0.1in}

\noindent $(\Leftarrow)$ If $\pi$ satisfies ~\eqref{eq:per cond} then $T^{(m+1)}(\pi)=\pi$ and therefore $\pi$ is periodic.

\end{proof}

\section{Counting the number of orbits of ${\mathcal P}_n$}

In~\cite{JB} it was shown how to apply Polya enumeration to count the number of distinct orbits of $T$ for a fixed value of $n$. We give a short exposition of this calculation. The idea is to identify each orbit of ${\mathcal P}_n$ with a necklace of black and white beads and then invoke Polya enumeration. 

If $n$ is triangular, then there is exactly one periodic point. So we assume that $n$ is not triangular. Therefore we assume that for some $m \in \N$, $n$ satisfies 
the inequality 
$$\frac{m(m+1)}{2} < n < \frac{(m+1)(m+2)}{2},$$
and we set $l = (m+1)(m+2)/2-n.$ 

Let $\pi$ be a periodic point of ${\mathcal P}_n$. By Theorem~\ref{Main-Result} we can identify $\pi$ with the arrangement of checkers and empty squares
on the $(m+1)$-th diagonal of the Ferrers graph. Furthermore by thinking of each empty square as being a white bead and each checker as a black bead we can identify $\pi$ with an arrangement of $l$ black beads and $(m+1-l)$ white beads. For example we  identify the periodic point $(5,4,4,2,2,1)$ with the following arrangement of black beads ($\bullet$) and white beads ($\circ$):
\begin{equation*}
(5,4,4,2,2,1) \,\, \leftrightarrow \,\,
\begin{array}{ccccccccccc}
\bullet & \bullet & \circ & \bullet & \circ & \circ 
\end{array}
\end{equation*}

Furthermore, we can view the action of $T$ on $\pi$ as rotating these black and white beads. Thus we have the cyclic group of order $(m+1)$, 
$C_{m+1}$, acting upon a line of $l$ black beads and $(m+1-l)$ white beads, and consequently we can identify the orbit, ${\mathcal O}_T(\pi)$, with a 
necklace consisting of $l$ black beads and $(m+1-l)$ white beads. For example, in the case of the periodic point $(5,4,4,2,2,1)$ we have the identification 
\begin{equation*}
{\mathcal O}((5,4,4,2,2,1)) \,\, \leftrightarrow \,\,
\left. \begin{array}{ccc}
 & \circ & \\ 
\bullet &  & \circ \\
\circ & & \bullet \\
 & \bullet & 
\end{array} \right.
\end{equation*}
These observations lead to the following theorem.

\begin{theorem}
The number of distinct orbits of ${\mathcal P}_n$ under the action of $T$ is equal to the number of necklaces consisting of $l$ black beads and 
$(m+1-l)$ white beads, where the symmetry group is the cyclic group of order $(m+1)$.
\end{theorem}

Counting the number of necklaces consisting of beads of two distinct colors, where the symmetry group is cyclic, is a standard exercise in Polya-enumeration. The interested reader is referred to~\cite[Chapter 27]{B} for the details. We simply state the basic result.

\begin{cor}
\label{No-Orbits}
Let $\varphi(d)$ denote the Euler phi function. The number of distinct orbits of ${\mathcal P}_n$ under the action of $T$ is the coefficient of the $b^lw^{m+1-l}$ term of the bivariate polynomial 
\begin{equation}\label{eq:count}
 \frac{1}{m+1}\sum_{d|(m+1)}\varphi(d)\left(b^d+w^d\right)^{(m+1)/d}.
\end{equation} 
\end{cor}

\section{A Slight Generalization of Bulgarian Solitaire}

Let $s=(s_n)$ be a sequence with $s_n \in \Z$. We define a map $T_s: {\mathcal P} \rightarrow {\mathcal P}$ in the following way using as a model our 
description of the $T$, see Observation~\ref{Act-T}.

\begin{proc}[Action of $T_s$] We start with a partition whose Ferrers graph is arranged on an infinite checkerboard. For each checker on $diag[k]$ we move the checker $(s_k \mod k)$ squares down the diagonal, looping around to the top of $diag[k]$ as needed.  Once we have moved all of the checkers diagonally, we check to see if there are any rows of checkers with spaces between checkers. For each such row we move the checkers to the left until all of the checkers are contiguous starting from the left side of the board.  We now check to see if there are columns of checkers with spaces between checkers. If there are we move the checkers up until we have a contiguous set of checkers starting at the top of the board. We now have the Ferrers graph of a partition, possibly the one with which we started. We note that our map $T$ (for Bulgarian Solitaire) is simply the special case $T_{(1,1,1,\ldots)}$.
\end{proc}

Our proof of Bulgarian solitaire immediately generalizes to give us the following result. The condition that $\gcd(k,s_k)=1$ is needed so that we can apply Lemma~\ref{CRTApp}.

\begin{theorem}[Generalized Bulgarian Solitaire]
\label{Gen-Result}
Let $\pi$ be a partition of $n$, and let $m$ be the integer such that $m(m+1)/2 \le n < (m+1)(m+2)/2$. Let $s=(s_k)$ be a sequence, with $s_n \in \Z$, satisfying the condition that $\gcd \left(k,s_k\right)=1$. Then $\pi$ is a periodic point of $T_s$ if and only if 
\begin{equation}
\# \textrm{ of checkers on } diag[k] = \left\{ \begin{array}{ll} k, & k=1,2, \ldots, m \\
n-m(m+1)/2, & k = m+1 \\
0, & k > m+1. \end{array}\right.
\end{equation} 
\end{theorem}

It is an interesting question to determine the periodic points and fixed points for $T_s$ when $s= \{s_n\}, s_n\in \Z$, is a sequence with   $\gcd(k,s_k) \not=1 $ for some $k$. In this case one can have periodic points (and fixed points) that are not periodic points (or fixed points) for 
Generalized Bulgarian Solitaire.  The trivial example is the sequence  $s=(1,2,3,4, \ldots)$. In this case every partition is a fixed point. A more interesting example is when $s$ is any integer valued sequence satisfying the condition $s_{15}=5$ and $s_{16}=8$. In this case the partition
$$ (15,14,13,11,11,10,9,9,6,6,5,4,3,1,1,1)$$ is a fixed point for $T_s$. This can be seen by looking at the last two rows of its rotated Ferrers graph.

\begin{equation*}
\begin{array}{c}
\bullet \circ \bullet \bullet \bullet \bullet \circ \bullet \bullet \bullet \bullet \circ \bullet \bullet \bullet  \\
 \bullet \circ  \circ \circ  \circ \circ \circ \circ \bullet  \circ \circ \circ \circ \circ \circ \circ \\ 
\end{array}
\end{equation*}


\begin{thebibliography}{99}

\bibitem{AE} G.~E.~Andrews and K.~Eriksson, {\it Integer Partitions}, Cambridge, 2004.

\bibitem{B} N.~L.~Biggs, {\it Discrete Mathematics}, 2nd ed., Oxford, 2002.

\bibitem{D} C.~Dickens, {\it Little Dorrit}, Everyman's Library, 1965.

\bibitem{H}B.~Hopkins, 30 years of Bulgarian Solitaire, to appear in {\it Coll. Math. Jour.}, Martin Gardner memorial issue.

\bibitem{JB} J.~Brandt, Cycles of Partitions, {\it Proc. Amer. Math. Soc.}, \textbf{85} (1982), 483-486.

\bibitem{MG} M.~Gardner, Mathematical Games: Tasks you cannot help finishing no matter how hard you try to block finishing them, {\it Scientific American}, 249 (1983), 12-21. Also available as Bulgarian Solitaire and Other Seemingly Endless Tasks, 27-43, in \emph{The Last Recreations}, Springer-Verlag, 2007

\bibitem{K1} S.~Limanov and A.~Toom, {\it Kvant}, (1980), No. 11, M655, 19-20.

\bibitem{S} A.~Sills, A Combinatorial Identity of A Partition Identity of Andrews and Stanley, {\it Int. J. of Math. and Math. Sci.} (2004), 2495-2503.

\bibitem{K2} A.~Toom, {\it Kvant} (1981), No. 7, 28-30.

\bibitem{VGRT} N.~B.~Vasilyev, V.~L.~Gutenmacher, Zh.~M.Rabbot, and A.~L.~Toom, {\it Mathematical Olympiads by Correspondence} (in Russian), Moscow Nauka, 1987, 144-147.

\end{thebibliography}
\end{document}